\documentclass[11pt]{article}

\usepackage[utf8]{inputenc}
\usepackage[T1]{fontenc}

\usepackage{amsmath, amssymb, amsthm, mathtools}

\usepackage[a4paper, margin=1in]{geometry}

\usepackage[colorlinks=true, linkcolor=blue, citecolor=blue]{hyperref}

\newtheorem{theorem}{Theorem}[section]
\newtheorem{proposition}[theorem]{Proposition}

\theoremstyle{remark}
\newtheorem{remark}[theorem]{Remark}

\title{
A new equivalence to the Riemann Hypothesis by means of the Salem integral equation
}

\author{
B.~J. Gonz\'{a}lez$^b$ \thanks{Email: bjglez@ull.edu.es} \and
E.~R. Negr\'{i}n$^{a,b}$ \thanks{Email: enegrin@ull.edu.es} \\
\small $^a$Departamento de An\'alisis Matem\'atico, Universidad de La Laguna. Spain \\
\small $^b$Instituto de Matem\'aticas y Aplicaciones (IMAULL), Universidad de La Laguna. Spain
}

\date{}

\begin{document}

\maketitle

% ---------- Abstract ----------
\begin{abstract}
This note presents a new equivalence to the Riemann Hypothesis by means of the Salem integral equation.
\end{abstract}

\bigskip


\medskip


% ============================================================
\section{Introduction}

\hspace{1 cm} The Riemann Hypothesis states that all non-trivial zeros on $0<\operatorname{Re} s<1$ of the Riemann zeta-funcion $\xi (s)$ lie in the line $\operatorname{Re} s=\frac{1}{2}$.

The Rapha\"el Salem publication of 1953 \cite{Salem} proves that the Riemann Hypothesis is true if and only if for each $\frac{1}{2}< \delta <1$ the only bounded measurable function $f(t)$ on $(0,\infty )$ satisfying
\[
\int_{0}^{\infty}\frac{t^{\delta -1}}{e^{xt}+1} f(t) dt =0,\quad \textrm{for\;all} \;x>0 ,
\] 
is $f=0$ almost everywhere.

In this note we prove that the Riemann Hypothesis is true if and only if for each $\frac{1}{2}<\delta <1$  the integral equation
\[
\int_{0}^{\infty}\frac{t^{\delta -1}}{e^{xt}+1} f(t) dt =0,\quad \textrm{for\;all}\;x>0,
\]
has not solutions of type $f(t) = t^{i\gamma}$, $t\in (0,\infty )$, $\gamma \in\mathbb{R}$.

Observe that $f(t) = t^{i\gamma}=e^{i\gamma \, \ln \, t}$ is a bounded measurable function on $(0,\infty )$ whose modulus is equal to one.

% ============================================================
%\section{Main Result}
\section*{Proof of the new equivalence}

\hspace{1 cm} For $\operatorname{Re} s>1$, the Riemann zeta-function $\xi (s)$  is represented by the convergent series:
\[
\xi (s)=\sum_{n=1}^{\infty}\frac{1}{n^s} .
\]

Also, denote by $\eta (s)$ the Dirichlet eta-function represented for $\operatorname{Re} s>0$ by the convergent series:
\[
\eta (s)= \sum_{n=1}^{\infty}\frac{(-1)^{n-1}}{n^s} .
\]

For $\operatorname{Re} s >1$ the next known relation holds
\[
\eta (s)=\left( 1- 2^{1-s}\right) \xi (s) ,
\]
which extends $\xi (s)$ to $0<\operatorname{Re} s<1$.

On the other hand, it is known that for $\operatorname{Re} s >0$
\begin{equation}\label{1.1}
\int_{0}^{\infty}\frac{t^{s -1}}{e^{xt}+1}  dt =x^{-s}\Gamma (s) \eta (s) ,\quad \textrm{for\;all} \: x>0.
\end{equation}

Thus, we obtain that $s= \delta + i\gamma$, $0<\delta < 1$, $\gamma\in\mathbb{R}$, is a zero of $\xi (s) $ if and only if 
\[
\int_{0}^{\infty}\frac{t^{\delta -1}}{e^{xt}+1} t^{i\gamma }  dt = 0,\quad \textrm{for\;all}\: x>0.
\]

Therefore, the Riemann Hypothesis is true if and only if for each $0<\delta < 1$, $\delta \neq \frac{1}{2}$, the integral equation
\[
\int_{0}^{\infty}\frac{t^{\delta -1}}{e^{xt}+1} f(t)  dt = 0, \quad \textrm{for\;all}\: x>0,
\]
has not solutions of type $f(t)=t^{i\gamma}$, $t\in (0,\infty )$, $\gamma\in \mathbb{R}$.

So, from the symmetry of the zeros of $\xi (s)$ on $0<\operatorname{Re} s<1$ with respect to the line $\operatorname{Re}  s =\frac{1}{2}$, one arrives to the next equivalence:

The Riemann Hypothesis is true if and only if for each $\frac{1}{2} < \delta < 1$ the integral equation 
\[
\int_{0}^{\infty}\frac{t^{\delta -1}}{e^{xt}+1} f(t)  dt = 0, \quad \textrm{for\;all}\: x>0,
\]
has not solutions of type $f(t)=t^{i\gamma}$, $t \in (0,\infty )$, $\gamma\in \mathbb{R}$.

\section*{A consequence}

\hspace{1 cm} Observe that from the Salem equivalence \cite{Salem} and from this new equivalence one has:

For each $0< \delta <1$, $\delta \neq \frac{1}{2}$, the only bounded measurable function $f(t)$ on $(0,\infty )$ satisfying
\[
\int_{0}^{\infty}\frac{t^{\delta -1}}{e^{xt}+1} f(t) dt =0,\quad \textrm{for\;all} \;x>0 ,
\] 
is $f =0$ almost everywhere if and only if for each $0<\delta <1$, $\delta \neq \frac{1}{2}$,
the integral equation
\[
\int_{0}^{\infty}\frac{t^{\delta -1}}{e^{xt}+1} f(t) dt =0,\quad \textrm{for\;all}\;x>0,
\]
has not solutions of type $f(t) = t^{i\gamma}$, $t\in (0,\infty )$, $\gamma \in\mathbb{R}$.

% ============================================================

\end{document}